\documentclass{article}
\usepackage{graphicx}
\begin{document}
\setlength{\textwidth}{5.75in}

\makeatletter
\def\@citex[#1]#2{\if@filesw\immediate\write\@auxout{\string\citation{#2}}\fi
  \def\@citea{}\@cite{\@for\@citeb:=#2\do
    {\@citea\def\@citea{, }\@ifundefined
       {b@\@citeb}{{\bf ?}\@warning
       {Citation `\@citeb' on page \thepage \space undefined}}%
\hbox{\csname b@\@citeb\endcsname}}}{#1}}
\makeatother


\makeatletter
\def\doublespaced{\baselineskip=\normalbaselineskip	
    \multiply\baselineskip by 150			
    \divide\baselineskip by 100}			
\def\doublespace{\doublespaced}				
\makeatother

\makeatletter
\def\mitspaced{\baselineskip=\normalbaselineskip	
    \multiply\baselineskip by 115			
    \divide\baselineskip by 100}			
\def\mitspace{\mitspaced}
\makeatother

\makeatletter
\def\capamispaced{\baselineskip=\normalbaselineskip	
    \multiply\baselineskip by 105			
    \divide\baselineskip by 100}			
\def\capamispace{\capamispaced}
\makeatother

\makeatletter
\def\singlespaced{\baselineskip=\normalbaselineskip}	
\def\singlespace{\singlespaced}				
\makeatother

\makeatletter
\def\triplespaced{\baselineskip=\normalbaselineskip	
    \multiply\baselineskip by 3}			
\makeatother

\makeatletter
\def\widenspacing{\multiply\baselineskip by 125		
    \divide\baselineskip by 100}			
\def\whitespace{\widenspacing}				
\makeatother

\font\eightmsam=msam8
\font\ninemsam=msam9
\font\tenmsam=msam10
\newtheorem{Theorem}{Theorem}[section]
\newtheorem{Corollary}[Theorem]{Corollary}
\newtheorem{Proposition}[Theorem]{Proposition}
\newtheorem{Lemma}[Theorem]{Lemma}
\newtheorem{Claim}[Theorem]{Claim}
\newtheorem{Definition}[Theorem]{Definition}
\newtheorem{Example}[Theorem]{Example}
\newcommand{\BB}{{\ninemsam\char'04}}
\renewcommand{\Box}{\mbox{\BB}\medskip}
\newcommand{\optimal}{$\Theta(n^{1/2})$}
\newcommand{\joptimal}{$\Theta(n^{1/j})$}
\newcommand{\nfourth}{$\Theta(n^{1/4})$}
\newcommand{\nsixth}{$\Theta(n^{1/6})$}
\newcommand{\nfourthlog}{$\Theta(n^{1/4}\log(n))$}
\newcommand{\linear}{$\Theta(n)$}
\newcommand{\logn}{$\Theta(\log(n))$}
\setlength{\baselineskip}{0.3in}
\hoffset -0.4in
\singlespace
\title{Digital Fixed Points, Approximate Fixed Points, and Universal Functions
}
\author{Laurence Boxer
         \thanks{
    Department of Computer and Information Sciences,
    Niagara University,
    Niagara University, NY 14109, USA;
    and Department of Computer Science and Engineering,
    State University of New York at Buffalo.
    E-mail: boxer@niagara.edu
    }
\and{
    Ozgur Ege
    \thanks{
    Department of Mathematics,
    Celal Bayar University,
    Muradiye,
    Manisa 45140, Turkey.
    E-mail: ozgur.ege@cbu.edu.tr
    }}
\and{
    Ismet Karaca
    \thanks{
    Department of Mathematics,
    Ege University,
    Bornova, Izmir 35100,
    Turkey.
    E-mail: ismet.karaca@ege.edu.tr
    }}
\and{
     Jonathan Lopez
     \thanks{
     Department of Mathematics,
     Canisius College,
     Buffalo, NY 14208, USA. 
     E-mail: lopez11@canisius.edu
     }}
\and{
Joel Louwsma
\thanks{
     Department of Mathematics,
     Niagara University, NY 14305, USA.
     E-mail: jlouwsma@niagara.edu
       }
    }
}
\date{ }
\maketitle

\begin{abstract}
A. Rosenfeld ~\cite{Rosenfeld}
introduced the notion of a digitally continuous
function between digital images, and showed that although digital images
need not have fixed point properties analogous to those of the Euclidean spaces
modeled by the images, there often are approximate fixed point properties
of such images. In the current paper, we obtain additional
results concerning fixed points and approximate fixed points of
digitally continuous functions. Among these are
several results concerning the relationship between
universal functions and the approximate fixed
point property (AFPP).


\end{abstract}

\section{Introduction}
In {\em digital topology}, we study geometric and topological
properties of digital images via tools adapted from geometric and algebraic
topology. Prominent among these tools is a digital version of
continuous functions. In the current paper, we study fixed points and
approximate fixed points of digitally continuous functions. 
We present a number of original results, as well as 
corrections of assertions that have appeared in previous
papers.

The paper is organized as follows.  Section~\ref{prelim}
reviews background material. In section~\ref{FPP}, we show
that a digital image $X$ has the Fixed Point Property (FPP)
if and only if $X$ has a single point. This result becomes
a key to many of the corrections we demonstrate in
section~\ref{correct-sec}. In section~\ref{approx-sec} we
introduce approximate fixed points and the Approximate
Fixed Point Property (AFPP). We give examples of digital images
that have, and that don't have, this property. In 
section~\ref{univ-funct} we study universal functions on
digital images and their relation to the AFPP. Concluding
remarks and an acknowledgement appear in 
sections~\ref{summary-sec} and~\ref{ack-sec}, respectively.

\section{Preliminaries}
\label{prelim}
\subsection{General Properties}
\label{dig-con}
A {\em fixed point} of a function $f: X \rightarrow X$ is a point
$x \in X$ such that $f(x)=x$.

For a finite set $X$, we denote by $|X|$ the number of distinct
members of $X$.

Let ${\bf N}$ be the set of natural numbers and let
${\bf Z}$ denote the set of integers.
Then ${\bf Z}^n$ is the set of lattice points in
Euclidean $n-$dimensional space.

A {\em digital image} is a pair $(X,\kappa)$, where
$\emptyset \neq X \subset {\bf Z}^n$ for some 
positive integer $n$ and
$\kappa$ is an adjacency relation on $X$. Technically, then,
a digital image $(X,\kappa)$ is an undirected graph whose
vertex set is the set of members of $X$ and whose edge set
is the set of unordered pairs $\{x_0,x_1\} \subset X$ such that
$x_0 \neq x_1$ and $x_0$ and $x_1$ are $\kappa$-adjacent.

Adjacency relations commonly used for
digital images include the following ~\cite{Kong}.
Two points $p$ and $q$ in ${\bf Z}^2$ are $8-adjacent$ if they
are distinct and differ by at most $1$ in each coordinate; $p$ and
$q$ in ${\bf Z}^2$ are $4-adjacent$ if they are 8-adjacent and
differ in exactly one coordinate.  Two points $p$ and $q$ in ${\bf Z}^3$
are $26-adjacent$ if they
are distinct and differ by at most $1$ in each coordinate; they are
$18-adjacent$ if they are 26-adjacent and
differ in at most two coordinates; they are
$6-adjacent$ if they are 18-adjacent and
differ in exactly one coordinate.  For $k \in \{4,8,6,18,26\}$,
a $k-neighbor$ of a lattice point $p$ is a point that is $k-$adjacent
to $p$.

The adjacencies discussed above are generalized as follows.
Let $u, n$ be positive integers, $1 \leq u \leq n$.
Distinct points $p,q \in {\bf Z}^n$ are called $c_u$-adjacent if
there are at most $u$ distinct coordinates $j$ for which
$|p_j-q_j| \, = \, 1$, and for all other coordinates $j$, $p_j = q_j$.
The notation $c_u$ represents the number of points $q \in {\bf Z}^n$
that are adjacent to a given point $p \in {\bf Z}^n$ in this sense. Thus the
values mentioned above:
if $n=1$ we have $c_1=2$; if $n=2$ we have $c_1=4$ and $c_2=8$;
if $n=3$ we have $c_1=6$, $c_2=18$, and $c_3=26$.
Yet more general adjacency relations are discussed in~\cite{Herman}.

Let $\kappa$ be an adjacency relation defined on ${\bf Z}^n$.
A digital image $X \subset {\bf Z}^n$ is $\kappa-connected$~\cite{Herman}
if and only if for
every pair of points $\{x,y\} \subset X$, $x \neq y$, there exists a set
$\{x_0, x_1,\ldots, x_c\} \subset X$ such that $x = x_0$, $x_c = y$,
and $x_i$ and $x_{i+1}$ are $\kappa-$neighbors, $i \in \{0,1,\ldots, c-1\}$.
A {\em $\kappa$-component} of $X$ is a maximal $\kappa$-connected subset of
$X$.

Often, we must assume some adjacency relation for
the {\em white pixels} in ${\bf Z}^n$,
{\em i.e.}, the pixels of ${\bf Z}^n \setminus X$ (the pixels that
belong to $X$ are sometimes referred to as {\em black pixels}).
In this paper, we are not concerned with adjacencies between
white pixels.

\begin{Definition} {\rm \cite{Boxer}}
Let $a, b \in {\bf Z}$, $a < b$.  A {\rm digital interval}
is a set of the form
\[ [a,b]_{{\bf Z}} ~=~ \{z \in {\bf Z} ~|~ a \leq z \leq b\} \]
in which $2-$adjacency is assumed. $\Box$
\end{Definition}

\begin{Definition}
\label{dig-cont}
{\rm (\cite{Boxer99}; see also ~\cite{Rosenfeld})}
Let $X \subset {\bf Z}^{n_0}$, $Y \subset {\bf Z}^{n_1}$.
Let $f:X\rightarrow Y$ be a function.
Let $\kappa_i$ be an adjacency relation defined on ${\bf Z}^{n_i}$,
$i \in \{0,1\}$.  We say $f$ is
{\rm $(\kappa_0,\kappa_1)-$continuous} if for every
$\kappa_0-$connected subset $A$ of $X$, $f(A)$ is a
$\kappa_1-$connected subset of $Y$. $\Box$
\end{Definition}

See also~\cite{Chen1,Chen2}, where similar
notions are referred to as {\em 
immersions}, {\em gradually varied operators},
and {\em gradually varied mappings}.

If $a$ and $b$ are members of a digital image $(X,\kappa)$,
we write $a \leftrightarrow_{\kappa} b$, or 
$a \leftrightarrow b$ when $\kappa$ is understood, to indicate
that either $a=b$ or $a$ and $b$ are $\kappa$-adjacent.

We say a function satisfying Definition~\ref{dig-cont}
is {\em digitally continuous}.
This definition implies the following.

\begin{Proposition}
\label{cont-connect}
{\rm (\cite{Boxer99}; see also~\cite{Rosenfeld})}
Let $X$ and $Y$ be digital images.  Then the function
$f: X \rightarrow Y$ is $(\kappa_0,\kappa_1)$-continuous if and only if
for every $\{x_0, x_1\} \subset X$ such that $x_0$ and $x_1$ are
$\kappa_0-$adjacent, $f(x_0) \leftrightarrow_{\kappa_1} f(x_1)$. 
$\Box$
\end{Proposition}

For example, if $\kappa$ is an adjacency relation on a digital image $Y$,
then $f:[a,b]_{{\bf Z}} \rightarrow Y$ is $(2,\kappa)-$continuous
if and only if for every
$\{c,c+1\} \subset [a,b]_{{\bf Z}}$, 
$f(c) \leftrightarrow_{\kappa} f(c+1)$.

We have the following.

\begin{Proposition}
\label{composition}
{\rm \cite{Boxer99}}
Composition preserves digital continuity, {\em i.e.},
if $f:X \rightarrow Y$ and $g:Y \rightarrow Z$ are, respectively,
$(\kappa_0,\kappa_1)-$continuous and $(\kappa_1,\kappa_2)-$continuous
functions, then the composite function $g \circ f: X \rightarrow Z$
is $(\kappa_0,\kappa_2)-$continuous. $\Box$
\end{Proposition}

We say digital images $(X,\kappa)$ and $(Y,\lambda)$ are
$(\kappa, \lambda)-isomorphic$ (called
$(\kappa, \lambda)-homeomorphic$ in ~\cite{Boxer,Boxer05}) if there
is a bijection $h: X \rightarrow Y$ that is $(\kappa,\lambda)$-continuous,
such that the function $h^{-1}: Y \rightarrow X$ is
$(\lambda, \kappa)$-continuous.

Classical notions of topology~\cite{Borsuk} 
yielded the concept of digital retraction
in~\cite{Boxer}.
Let $(X, \kappa)$ be a digital image and let $A$ be a nonempty subset of $X$.
A {\em retraction} of $X$ onto $A$ is a $(\kappa,\kappa)$-continuous function
$r:X \rightarrow A$ such that $r(a)=a$ for all $a \in A$. 

A {\em digital simple closed curve} is a digital image
$X=\{x_i\}_{i=0}^{m-1}$, with $m \geq 4$, such that
the points of $X$ are {\em labeled circularly}, i.e.,
$x_i$ and $x_j$ are adjacent
if and only if $j=(i-1) ~(mod ~m)$ or $j=(i+1) ~(mod ~m)$.

\subsection{Digital homotopy}
A homotopy between continuous functions may be thought of as
a continuous deformation of one of the functions into the other
over a time period.

\begin{Definition}{\rm (\cite{Boxer99}; see also \cite{Khalimsky})}
\label{htpy-2nd-def}
Let $X$ and $Y$ be digital images.
Let $f,g: X \rightarrow Y$ be $(\kappa,\kappa')$-continuous functions.
Suppose there is a positive integer $m$ and a function
$F: X \times [0,m]_{{\bf Z}} \rightarrow Y$
such that
\begin{itemize}
\item for all $x \in X$, $F(x,0) = f(x)$ and $F(x,m) = g(x)$;
\item for all $x \in X$, the induced function
      $F_x: [0,m]_{{\bf Z}} \rightarrow Y$ defined by
          \[ F_x(t) ~=~ F(x,t) \mbox{ for all } t \in [0,m]_{{\bf Z}} \]
          is $(2,\kappa')-$continuous.
\item for all $t \in [0,m]_{{\bf Z}}$, the induced function
         $F_t: X \rightarrow Y$ defined by
          \[ F_t(x) ~=~ F(x,t) \mbox{ for all } x \in  X \]
          is $(\kappa,\kappa')-$continuous.
\end{itemize}
Then $F$ is a {\rm digital $(\kappa,\kappa')-$homotopy between} $f$ and
$g$, and $f$ and $g$ are {\rm digitally $(\kappa,\kappa')-$homotopic in} $Y$.
$\Box$
\end{Definition}

When the adjacency relations $\kappa$ and $\kappa'$ are understood in context,
we say $f$ and $g$ are {\em digitally homotopic} to abbreviate
``digitally $(\kappa,\kappa')-$homotopic in $Y$."






\begin{Definition}
\label{htpy-trivial}
A digital image $(X,\kappa)$ is
$\kappa$-{\rm contractible \cite{Khalimsky,Boxer}}
if its identity map is $(\kappa, \kappa)$-homotopic to a constant
function $\overline{p}$ for some $p \in X$.
$\Box$
\end{Definition}

When $\kappa$ is understood, we speak of {\em contractibility}
for short.



\subsection{Digital simplicial homology}
Our presentation of digital simplicial homology is taken
from that of~\cite{Ege&Karaca3}.

A set of $m+1$ distinct mutually adjacent points
is an {\em $m$-simplex}.

\begin{Definition}
\label{alpha-def}
If $\alpha_q$ is the number of $(\kappa,q)$-simplices in $X$
and $m = \max\{q \in {\bf N}^* \, | \, \alpha_q \neq 0\}$, then
$m$ is the {\em dimension of} $(X,\kappa)$, denoted
$dim(X, \kappa)$ or $dim(X)$, and 
the {\em Euler characteristic of $(X,\kappa)$},
$\chi(X,\kappa)$, is defined by
\[ \chi(X,\kappa)= \sum_{q=0}^m (-1)^q \alpha_q.~ \Box \]
\end{Definition}

For $q \in {\bf N}^*$, the group of $q$-chains of
$(X,\kappa)$, denoted $C_q^{\kappa}(X)$, is the free
Abelian group with basis being the set of $q$-simplices of
$X$.

Let 
$\delta_q: C_q^{\kappa}(X) \rightarrow C_{q-1}^{\kappa}(X)$ 
defined by
\[ \delta_q(<p_0,p_1,\ldots,p_q>)= \left \{ \begin{array}{ll}
   \sum_{i=0}^q (-1)^i <p_0,p_1,\ldots, \hat{p}_i, \ldots,p_q>
   & \mbox{if } 0 \leq q \leq dim(X); \\
   0 & \mbox{if } q > dim(X),
   \end{array} \right .
\]
where $\hat{p}_i$ means that $p_i$ is omitted from the
vertices of the simplex considered. Then $\delta_q$ is a
homomorphism, and we have
$\delta_{q-1} \circ \delta_q = 0$~\cite{AKO}. This gives
rise to the following groups~\cite{BKO}.
\begin{itemize}
\item $Z_q^{\kappa}(X)=Ker \, \delta_q$, the group of digital 
       simplicial $q$-{\em cycles} of $X$.
\item $B_q^{\kappa}(X)=Im \, \delta_{q+1}$, the group of    
       digital simplicial $q$-{\em boundaries} of $X$.
\item The quotient group 
      $H_q^{\kappa}(X)=Z_q^{\kappa}(X) \, / \, B_q^{\kappa}(X)$,
      the $q$-th digital simplicial homology group of $X$.
\end{itemize}

We have the following.

\begin{Theorem}
\label{homol-props}
{\rm\cite{BKO}}
Let $(X, \kappa)$ be a directed digital simplicial 
complex of dimension $m$.
\begin{itemize}
\item $H_q(X)$ is a finitely generated abelian group for 
      every $q \geq 0$.
\item $H_q(X)$ is a trivial group for all $q > m$.
\item $H_m(X)$ is a free abelian group, possibly $\{0\}$. $\Box$
\end{itemize}
\end{Theorem}

\section{Fixed point property}
\label{FPP}
We say a digital image $(X,\kappa)$ has the
{\em fixed point property (FPP)} if every
$(\kappa,\kappa)$-continuous function $f: X \rightarrow X$ has a
fixed point. Some properties of digital images with the FPP were
studied in~\cite{Ege&Karaca}. However, the following shows that
for digital images with $c_u$-adjacencies,
the FPP is not very interesting.

\begin{Theorem}
\label{FPP-triv}
Let $(X,\kappa)$ be a digital image. Then
$(X, \kappa)$ has the FPP if and only if $|X|=1$.
\end{Theorem}

{\em Proof}: Clearly, if $|X|=1$ then $(X,\kappa)$ has the FPP.

Now suppose $|X|>1$.
If $(X,\kappa)$ has more than 1 $\kappa$-component, then
there is a $(\kappa,\kappa)$ continuous map $f: X \rightarrow X$
such that for all $x \in X$, $x$ and $f(x)$ are in different $\kappa$-components
of $X$. Such a map does not have a fixed point.

Therefore, we may assume $X$ is $\kappa$-connected. Since $|X| > 1$, there
are distinct $\kappa$-adjacent points $x_0,x_1 \in X$. 
Consider the map $f: X \rightarrow X$ given by
\[ f(x) = \left \{ \begin{array}{ll}
          x_0 & \mbox{if } x \neq x_0; \\
          x_1 & \mbox{if } x = x_0.
         \end{array} \right . \]
Consider a pair $y_0,y_1$ of $\kappa$-adjacent members of $X$.
\begin{itemize}
\item If one of these points, say, $y_0$, coincides with $x_0$, we have
      $f(y_0)=f(x_0)=x_1$ and, since $y_1 \neq x_0$, $f(y_1)=x_0$, so
      $f(y_0)$ and $f(y_1)$ are $\kappa$-adjacent.
\item If both $y_0$ and $y_1$ are distinct from $x_0$, then
      $f(y_0)=x_1=f(y_1)$
\end{itemize}
Therefore, $f$ is $(\kappa,\kappa)$-continuous. Clearly, $f$ does not have a
fixed point. Therefore, $(X,\kappa)$ does not have the FPP.
$\Box$

\section{Corrections of published assertions}
\label{correct-sec}
In this section, we correct some assertions that appear in~\cite{Ege&Karaca,Ege&Karaca3}.

We show below that the function 
$F: [0,1]_{\bf Z} \rightarrow [0,1]_{\bf Z}$ defined by
$F(x)=1-x$ (i.e., $F(0)=1$, $F(1)=0$) provides a 
counterexample to several of the assertions of
~\cite{Ege&Karaca}. Clearly this function is
$(2,2)$-continuous and does not have a fixed point.

We will need the following.

\begin{Definition}
{\rm \cite{Ege&Karaca2}}
Let $(X, \kappa)$ be a digital image whose digital homology
groups are finitely generated and vanish above some
dimension $n$. Let $f: X \rightarrow X$ be a
$(\kappa,\kappa)$-continuous map. The 
{\rm Lefschetz number of} $f$, denoted $\lambda(f)$, is
defined as
\[ \lambda(f)=\sum_{i=0}^n (-1)^i \, tr(f_{i,*}), \]
where $f_{i,*}: H_i^{\kappa}(X) \rightarrow H_i^{\kappa}(X)$
is the map induced by $f$ on the $i^{th}$ homology group of
$(X,\kappa)$ and $tr(f_{i,*})$ is the trace of $f_{i,*}$. $\Box$
\end{Definition}

In studying digital maps from a sphere to itself, there is
a question of how to represent a Euclidean sphere digitally.
\begin{itemize}
\item As in ~\cite{Ege&Karaca3}, we will represent $S^1$ by the
      set $S_1 = [-1,1]_{\bf Z}^2 \setminus \{(0,0)\} 
      \subset {\bf Z}^2$ and $c_1$-adjacency with points
      $\{x_j\}_{j=0}^7$ labeled circularly.
\item More generally, as in ~\cite{Ege&Karaca3}, we will represent $S^n$ by the
      set $S_n = [-1,1]_{\bf Z}^{n+1} \setminus \{0_{n+1}\} 
      \subset {\bf Z}^{n+1}$ and $c_1$-adjacency,
      where $0_{n+1}$ is the origin in
      ${\bf Z}^{n+1}$.
\end{itemize}

\begin{Definition}
\label{degree}
{\rm \cite{Ege&Karaca3,Spanier}}
Suppose a continuous function 
$f: (S_n,\kappa) \rightarrow (S_n,\kappa)$
induces a homomorphism on the $n$-th homology group,
$f_*: H_n^{\kappa}(S_n) \rightarrow H_n^{\kappa}(S_n)$, such that
$f_*([x])=m[x]$ for some fixed $m \in {\bf Z}$, 
where $[x]$ is a generator of $H_n^{\kappa}(S_n)$. 
The value of $m$ is the {\em degree of} $f$. $\Box$
\end{Definition}

\begin{Theorem}
\label{scc-iso-htpy}
{\rm ~\cite{Boxer10}} Let $S$ be a digital simple
closed curve. For an isomorphism of $S$ and a 
continuous non-surjective self-map of $S$ to be 
homotopic, we must have $|S|=4$. $\Box$
\end{Theorem}

We state the following corrections.
\begin{itemize}
\item Incorrect assertion stated as Theorem 3.3 of
      ~\cite{Ege&Karaca}:
      If $(X,\kappa)$ is a finite digital image and
      $f: X \rightarrow X$ is a $(\kappa,\kappa)$-continuous
      function with $\lambda(f) \neq 0$, then $f$ has a
      fixed point.
      
      In fact, the function $F$ defined above is a 
      counterexample to this assertion, since it is
      easily seen that $\lambda(F) \neq 0$.
\item Incorrect assertion stated as Theorem 3.4 of ~\cite{Ege&Karaca}:
      Every $(c_1,c_1)$-continuous function
      $f: [0,1]_Z \rightarrow [0,1]_Z$ has a fixed point.

      In fact, ~\cite{Rosenfeld} shows that this assertion
      is false, and the function $F$ defined above
      is a counterexample.
\item Incorrect assertion stated as Theorem 3.5 of ~\cite{Ege&Karaca}:
      Let $X = \{(0,0), (1,0), (0,1), (1,1)\} \subset Z^2$. Then every
      $(c_1,c_1)$-continuous function $f: X \rightarrow X$
      has a fixed point.

      In fact, we can use the function $F$ above to obtain a counterexample.
      Let $G: X \rightarrow X$ be defined by
      $G(x,y)=(x,F(y))$. Then $G$ is $(c_1,c_1)$-continuous and has no fixed point. Alternately, it follows from
      Theorem~\ref{FPP-triv} that the assertion is incorrect.
\item Incorrect assertion stated as Theorem 3.8 of    
      ~\cite{Ege&Karaca}:
      Let $(X,\kappa)$ be a
      $\kappa$-contractible digital image. Then every
      $(\kappa,\kappa)$-continuous map $f: X \rightarrow X$ has a fixed point.

      In fact, since $[0,1]_Z$ is $c_1$-contractible, the function $F$
      above provides a counterexample to this assertion.
      Alternately, it follows from
      Theorem~\ref{FPP-triv} that the assertion is incorrect.
\item Incorrect assertion stated as Example 3.9 of
      ~\cite{Ege&Karaca}:
      Let $X=\{(0,0), (0,1), (1,1)\} \subset {\bf Z}^2$.
      Then $(X,c_2)$ has the FPP.
      
      In fact, the map $f: X \rightarrow X$ defined by
      $f(0,0)=(0,1)$, $f(0,1)=(1,1)$, $f(1,1)=(0,0)$ is
      $(c_2,c_2)$-continuous and has no fixed points. 
      Alternately, it follows from
      Theorem~\ref{FPP-triv} that the assertion is incorrect.
\item Incorrect assertion stated as Corollary 3.10 of
      ~\cite{Ege&Karaca}: Any digital image with the same
      digital homology groups as a single point image has
      the FPP.
      
      To show this assertion is incorrect, observe that if
      $X=\{(0,0), (0,1), (1,1)\} \subset {\bf Z}^2$ and $Y$
      is a digital image of one point in ${\bf Z}^2$, then
      ~\cite{AKO,BKO}
      \[ H_q^8(X) = H_q^8(Y) = 
         \left \{ \begin{array}{ll}
                   {\bf Z} & \mbox{if } q=0; \\
                   0 & \mbox{if } q \neq 0.
                   \end{array} \right .                  
      \]
      It follows from Theorem~\ref{FPP-triv} that the
      assertion is incorrect.
\item Incorrect assertions stated as Example~3.17 and
      Corollary~3.18 of~\cite{Ege&Karaca}: The digital 
      images $(MSS_6',6)$ and $(P^2,6)$, each with more
      than one point, have the FPP.
      
      It follows from
      Theorem~\ref{FPP-triv} that these assertions are
      incorrect.
\item Incorrect assertion stated as Theorem~3.5 of
      ~\cite{Ege&Karaca3}:
      If $(X,\kappa)$ is a finite digital image and
      $f: X \rightarrow X$ is a $(\kappa,\kappa)$-continuous
      function with $\lambda(f) \neq 0$, then any map
      homotopic to $f$ has a fixed point.
      
      In fact, we observed above that the function $F$, which is 
      homotopic to itself and has $\lambda(F) \neq 0$, does not 
      have a fixed point.
 
\item Incorrect assertion stated as Theorem 3.7 of
      ~\cite{Ege&Karaca3}: If $(X,\kappa)$ is a digital image
      such that $\chi(X,\kappa) \neq 0$, then any map homotopic
      to the identity has a fixed point.
      
      In fact, we can take $X = [0,1]_{\bf Z}$, for which
      $\chi(X,c_1)=(-1)^1 (2) + (-1)^2 (1) \neq 0$, and the 
      function $F$ discussed above is homotopic to $1_X$ and
      does not have a fixed point.
\item Incorrect assertion stated as Theorem~3.11 of
      ~\cite{Ege&Karaca3}: Let $(S_n,c_1) \subset {\bf Z}^{n+1}$
      be a digital $n$-sphere as described above, where 
      $n \in \{1,2\}$. If $f: S_n \rightarrow S_n$ is a 
      continuous map of degree $m \neq 1$, then $f$ has 
      a fixed point.
      
      In fact, we have the following.
      Elementary calculations show that
      $H_1^{c_1}(S_1) \approx {\bf Z}$; also,
      $H_1^{c_1}(S_2) \approx {\bf Z}^{23}$~\cite{D&K}.
      For $n \in \{1,2\}$, as in the proof of Theorem 3.1, we 
      can choose distinct and adjacent $x_0$ and $x_1$ in $S_n$
      and let $f:S_n\rightarrow S_n$ be given by $f(x)=x_0$ for 
      $x\neq x_0$ and $f(x_0)=x_1$.  Clearly, $f$ is 
      continuous and does not have a fixed point.  Since 
      $f_*:H_1(S_n)\rightarrow H_1(S_n)$ is $0$, the degree 
      of $f$ is $0$.
\item Proposition~3.12 of ~\cite{Ege&Karaca3} depends
      on an unstated assumption that (recall
      Definition~\ref{alpha-def})
      $\alpha_q^{\kappa}(X)$ is finite for all $q$,
      a condition that is satisfied if and only if
      $X$ is finite;
      after all, one can study infinite digital images 
      $(X,\kappa)$, as in~\cite{Boxer07}, for which, 
      e.g., $\alpha_0^{\kappa}(X) = \infty$. E.g.,
      we could take $X={\bf Z}$; according to
      Definition~\ref{alpha-def},
      $\chi({\bf Z},c_1)$ is undefined, since
      $\alpha_1^{c_1}({\bf Z})=
      \alpha_0^{c_1}({\bf Z})=\infty$.
      Therefore, the proposition should be stated as
      follows.
      \begin{quote}
      {\em Let $(X,\kappa)$ be a finite digital image and 
      suppose $f: (X,\kappa) \rightarrow (X,\kappa)$ 
      is continuous.
      If $f_*: H_*^{\kappa}(X) \rightarrow H_*^{\kappa}(X)$
      is defined by $f_*(z)=kz$
      where $k \in {\bf Z}$,
      i.e., if there exists $k \in {\bf Z}$ such
      that in every dimension $i$ we have
      $f$ inducing the homomorphism
      $f_{i*}: H_i^{\kappa}(X) \rightarrow H_i^{\kappa}(X)$ 
      defined by $f_{i*}(z)=kz$,
      then $\lambda(f) = k \, \chi(X)$.}
      \end{quote}
\item A theoretically minor, but possibly
      confusing, error in Theorem~3.14 
      of~\cite{Ege&Karaca3}: In discussing an
      {\em antipodal map} $f: X \rightarrow X$, one 
      needs the property that
      for every $x \in X$ we have $-x \in X$; this
      property does not characterize the version of
      $S_2$ used in Theorem~3.14 
      of~\cite{Ege&Karaca3}. In the following,
      we use
      $S_2 = [-1,1]_{\bf Z}^3 \setminus \{(0,0,0)\}$, as described above.
      
      Theorem~3.14 of ~\cite{Ege&Karaca3} 
      asserts that
      \begin{quote}
      If $\alpha_i: (S_i,c_1 ) \rightarrow (S_i,c_1 )$ is the antipodal map between two
      digital $i$-spheres $S_i \subset Z^{i+1}$,
      for $i \in \{1,2\}$, then $\alpha_i$ has 
      degree $(-1)^(i+1)$.
      \end{quote}
      
      In fact, we show that this assertion is 
      correct for $i=1$, although an argument
      different from that of~\cite{Ege&Karaca3} 
      must be given, as the argument 
      of~\cite{Ege&Karaca3} makes use of 
      Theorem~3.4 of~\cite{Ege&Karaca3} 
      (= Theorem~3.3 of [14]), which, as noted 
      above, is incorrect. For $i=2$, we show the 
      assertion is not well defined.
      \begin{itemize}
      \item For $i=1$, we have the following.
      Let the points $\{e_j\}_{j=0}^7$ of $S_1$ be 
      circularly ordered. For notational 
      convenience, let $e_8=e_0$, and, more 
      generally, index arithmetic is assumed to be 
      modulo 8. The generators of the 1-chains of 
      $S_1$ are the members of 
      $\{<e_j e_{j+1}>\}_{j=0}^7$.  We have
      \[ 0=\delta(\sum_{j=0}^7 u_j <e_j e_{j+1}>) =
         \sum_{j=0}^7 u_j(e_j - e_{j+1})=
         \sum_{j=1}^8(u_j - u_{j-1})e_j
      \]
      implies $u_0=u_1=\cdots = u_7$. Therefore,
      $Z_1(S_1)$ is generated by
      $\sum_{j=0}^7 <e_j e_{j+1}>$. Since, clearly,
      $B_1(S_1)=\{0\}$, we have 
      $H_1(S_1) = Z_1(S_1)/ B_1(S_1)$ is isomorphic
      to ${\bf Z}$. Therefore, the homomorphism 
      $(\alpha_1)_*: H_1(S_1 )\rightarrow H_1(S_1)$  
      induced by $\alpha_1$ must satisfy 
      $(\alpha_1)_*(x)=kx$ for some integer $k$.
      
      Indeed, since the antipode of $e_j$ is
      $e_{j+4}$, we have
      \[ (\alpha_1)_*(\sum_{j=0}^7 <e_j e_{j+1}>) =
      \sum_{j=0}^7 <e_{j+4} e_{j+5}> =
      \sum_{j=0}^7 <e_j e_{j+1}>,
      \]
      so $k=1=(-1)^{1+1}$, as asserted.
      
      \item For i=2, we observe that, using the 
      $c_1$ adjacency, there is no triple of distinct, 
      mutually adjacent points in $S_2$. Therefore, 
      $H_2(S_2 )=\{0\}$. Therefore, the degree of 
      $\alpha_2$ is not well defined since for any 
      integer $k$ we have $(\alpha_2)_*(x)=kx$ for 
      all $x=0 \in H_2(S_2)$.
      \end{itemize}
\item Incorrect assertion stated as Theorem 3.15 of
      ~\cite{Ege&Karaca3}: Let $S_1$ be a digital
      simple closed curve in ${\bf Z}^2$. If 
      $h: (S_1,c_1) \rightarrow (S_1,c_1)$ is a continuous
      function that is homotopic to a constant function in
      $S_1$, then $h$ has a fixed point.
      
      In fact, we can take $S_1$ as above with its
      points ordered circularly, $S_1= \{x_j\}_{j=0}^7$ where
      distinct points $x_u,x_v$ are adjacent if and only if
      $u+1=v~mod~8$ or $u-1=v~mod~8$. Then, as in the proof
      of Theorem~\ref{FPP-triv}, the function
      $h: S_1 \rightarrow S_1$ given by
      \[ h(x)= \left \{ \begin{array}{ll}
      x_0 & \mbox{if } x \neq x_0; \\
      x_1 & \mbox{if } x=x_0,
      \end{array} \right . \]
      is continuous and homotopic to the constant function
      $\overline{x_0}$ in
      $S_1$ but has no fixed point.
\item Correct (for $|S_1|>4$) assertion incorrectly
      ``proven" as Corollary
      3.16 of ~\cite{Ege&Karaca3}: Let $S_1$ be as above. Let
      $h: (S_1,c_1) \rightarrow (S_1,c_1)$ be given by
      $h(x_i)=x_{(i+1)~mod~m}$, where $m=|S_1|$. Then
      $h$ is not homotopic in $S_1$ to a constant map.
      
      The argument given for this assertion 
      in~\cite{Ege&Karaca3} depends on Theorem 3.15 of
      ~\cite{Ege&Karaca3}, which, we have shown above, is
      incorrect.  However, since $h$ is easily seen to 
      be an isomorphism, by Theorem~\ref{scc-iso-htpy},
      the current assertion is true if and only if 
      $|S_1|>4$.     
\end{itemize}

\section{Approximate fixed points}
\label{approx-sec}
Given a digital image $(X,\kappa)$ and a $(\kappa,\kappa)$-continuous
function $f: X \rightarrow X$, we say $p \in X$ is an
{\em approximate fixed point of} $f$ if either $f(p)=p$, or $p$ and $f(p)$ are
$\kappa$-adjacent.
We say a digital image $(X,\kappa)$ has the
{\em approximate fixed point property (AFPP)} if every
$(\kappa,\kappa)$-continuous function $f: X \rightarrow X$ has an
approximate fixed point.

\begin{Theorem}
{\rm ~\cite{Rosenfeld}}
\label{AR-thm}
Let $I = \Pi_{i=1}^n \, [a_i,b_i]_{{\bf Z}}$.
Then $(I,c_n)$ has the AFPP.
$\Box$
\end{Theorem}

Theorems~\ref{FPP-triv} and~\ref{AR-thm} show that it is worthwhile to consider the AFPP, rather
than the FPP, for digital images.  We have the following.

\begin{Theorem}
\label{iso-afp}
Suppose $(X, \kappa)$ has the AFPP and there is a
$(\kappa, \lambda)$-isomorphism $h: X \rightarrow Y$.  Then
$(Y, \lambda)$ has the AFPP.
\end{Theorem}

{\em Proof}:  Let $f: Y \rightarrow Y$ be $(\lambda,\lambda)$-continuous.
By Proposition~\ref{composition}, the function
$g \, = \, h^{-1} \circ f \circ h: X \rightarrow X$
is $(\kappa, \kappa)$ continuous,
so our hypothesis implies there exists $p \in X$ such that 
$p \leftrightarrow_{\kappa} g(p)$.

Then
\[ h(p) \leftrightarrow_{\lambda} h \circ g(p) \, = \, h \circ h^{-1} \circ f \circ h(p) \,
    = \, f(h(p)), \]
so $h(p)$ is an approximate fixed point of $f$.
$\Box$

\begin{Proposition}
A digital simple closed curve of 4 or more points does not
have the AFPP.
\end{Proposition}

{\em Proof}: Let $S=(\{s_i\}_{i=0}^{m-1}, \kappa)$ with
$s_i$ and $s_j$ $\kappa$-adjacent if and only if
$i=(j+1) \, mod \, m$ or $i=(j-1) \, mod \, m$. Then the
function $f: S \rightarrow S$ defined by
$f(s_i)= s_{(i+2) \, mod \, m}$ is 
$(\kappa,\kappa)$-continuous, and, for each $i$, $s_i$ and
$f(s_i)$ are neither equal nor $\kappa$-adjacent. $\Box$

Next, we show retractions preserve the AFPP.

\begin{Theorem}
\label{approx-fp-thm}
Let $(X,\kappa)$ be a digital image, and let $Y \subset X$ be a
$(\kappa,\kappa)$-retract of $X$.  If $(X,\kappa)$ has the AFPP,
then $(Y,\kappa)$ has the AFPP.
\end{Theorem}

{\em Proof}:  Let $r: X \rightarrow Y$ be a $(\kappa, \kappa)$ retraction.
Let $f: Y \rightarrow Y$ be a $(\kappa, \kappa)$-continuous function.
Let $i: Y \rightarrow X$ be the inclusion map.
By Proposition~\ref{composition},
$g=i \circ f \circ r: X \rightarrow X$ is $(\kappa,\kappa)$-continuous.
Therefore, $g$ has an approximate fixed point $x_0 \in X$.

Let $x_1 = g(x_0) \in Y$. By choice of $x_0$,
      it follows that $x_0 \leftrightarrow x_1$.  Then
\[ x_1 ~=~ g(x_0) \leftrightarrow g(x_1) ~=~ i \circ f \circ   
   r(x_1) ~=~ i \circ f(x_1) ~=~ f(x_1).\]
Thus $x_1$ is an approximate fixed point of $f$.
      $\Box$

Following a classical construction of topology, the
{\em wedge} of two digital images $(A, \kappa)$ and
$(B, \lambda)$, denoted $A \wedge B$, is defined
~\cite{Han} as the union of the digital images
$(A',\mu)$ and $(B', \mu)$, where
\begin{itemize}
\item $A' \cap B'$ has a single point, $p$;
\item If $a \in A'$ and $b \in B'$ are $\mu$-adjacent, then
      either $a=p$ or $b=p$;
\item $(A', \mu)$ and $(A, \kappa)$ are isomorphic; and
\item $(B', \mu)$ and $(B,\lambda)$ are isomorphic.
\end{itemize}
In practice, we often have $\kappa=\lambda=\mu$,
$A=A'$, $B=B'$.

We have the following.

\begin{Theorem}
Let $A$ and $B$ be digital images.  Then
$(A \wedge B, \kappa)$ has the AFPP if and only if both
$(A,\kappa)$ and $(B, \kappa)$ have the AFPP.
\end{Theorem}

{\em Proof}: Let $A \cap B = \{p\}$. Let
$p_A, p_B: A \wedge B \rightarrow A \wedge B$ be the
functions
\[ p_A(x) = \left \{ \begin{array}{ll}
                      x & \mbox{if } x \in A; \\
                      p & \mbox{if } x \in B.
                     \end{array} \right .
            ~~~~~~~~
   p_B(x) = \left \{ \begin{array}{ll}
                      p & \mbox{if } x \in A; \\
                      x & \mbox{if } x \in B.
                     \end{array} \right .
\]
It is easily seen that both of these functions are well
defined and $(\kappa,\kappa)$-continuous. Also, let
$i_A: A \rightarrow A \wedge B$ and 
$i_B: B \rightarrow A \wedge B$ be the inclusion functions,
which are clearly $(\kappa,\kappa)$-continuous.

Suppose $(A,\kappa)$ and $(B,\kappa)$ have the AFPP. Let
$f: A \wedge B \rightarrow A \wedge B$ be 
$(\kappa,\kappa)$-continuous. We must show that there
exists a point of $A \wedge B$ that is equal or
$\kappa$-adjacent to its image under $f$. If $f(p)=p$, then
we have realized that goal. Otherwise, without loss of
generality, $f(p) \in A \setminus \{p\}$. By
Proposition~\ref{composition},
$h = p_A \circ f \circ i_A: A \rightarrow A$ is
$(\kappa, \kappa)$-continuous. Since $A$ has the AFPP,
there exists $a \in A$ such that
\begin{equation}
\label{approx-fp}
   h(a) \leftrightarrow_{\kappa} a.
\end{equation}
If $f(a) \in B$, then
\[ p=p_A \circ f(a)=p_A \circ f \circ i_A(a)=h(a). \]
It follows from statement~(\ref{approx-fp}) that 
\begin{equation}
\label{p&a}
p \leftrightarrow_{\kappa} a.
\end{equation}
If $f(a) \neq p$ then $f(a) \in B \setminus \{p\}$ and 
$f(p) \in A \setminus \{p\}$, so
$f(a)$ and $f(p)$ are distinct, non-adjacent points.
This is a contradiction of statement~(\ref{p&a}), since $f$
is continuous.
Therefore, we have $f(a) \in A$. Then
\[ f(a)=p_A \circ f(a) = p_A \circ f \circ i_A(a) = h(a).
\]
It follows from statement~(\ref{approx-fp}) that 
$f(a) \leftrightarrow_{\kappa} a$.

Since $f$ was arbitrarily selected, it follows that
$(A \wedge B, \kappa)$ has the AFPP.

Conversely, suppose $(A \wedge B, \kappa)$ has the 
AFPP. Since the maps $p_A$ and $p_B$ are 
$(\kappa, \kappa)$-retractions of
$(A \wedge B, \kappa)$ onto $(A,\kappa)$ and
$(B,\kappa)$, respectively, it follows from
Theorem~\ref{approx-fp-thm}
that $(A,\kappa)$ and
$(B,\kappa)$ have the AFPP. $\Box$

\section{Universal functions and the AFPP}
\label{univ-funct}
In this section, we define the notion of a universal
function and study its relation to the AFPP.

\begin{Definition}
\label{universal}
Let $(X,\kappa)$ and $(Y,\lambda)$ be digital images.
A $(\kappa,\lambda)$-continuous function
$f: X \rightarrow Y$ is {\rm universal for} $(X,Y)$
if given a $(\kappa,\lambda)$-continuous function
$g: X \rightarrow Y$, there exists $x \in X$ such that
$f(x) \leftrightarrow_{\lambda} g(x)$.
\end{Definition}

The notion of a dominating set in graph theory corresponds
to the notion of a dense set in a topological space.

\begin{Definition}
\label{dominating}
{\rm \cite{C&L}}
Let $(X,\kappa)$ be a nonempty digital image. Let $Y$
be a nonempty subset of $X$. We say $Y$ is
{\rm $\kappa$-dominating} in $X$ if for
every $x \in X$ there
exists $y \in Y$ such that 
$x \leftrightarrow_{\kappa} y$.
\end{Definition}

\begin{Theorem}
\label{univ-dominating}
Let $(X,\kappa)$ and $(Y,\lambda)$ be digital images.
Let $f: X \rightarrow Y$ be a universal function
for $(X,Y)$. Then $f(X)$ is $\lambda$-dominating in $Y$.
\end{Theorem}

{\em Proof}: Let $y \in Y$ and consider the constant
function $c_y: X \rightarrow Y$ defined by
$c_y(x)=y$ for all $x \in X$. This function is clearly
$(\kappa,\lambda)$ continuous. Since $f$ is universal,
there exists $x_y \in X$ such that $f(x_y)$ is either
equal to or $\lambda$-adjacent to $y$. Since $y$
was arbitrarily chosen, the assertion follows.
$\Box$

\begin{Proposition}
\label{non-dominating-example}
Let $X$ be a $\kappa$-connected digital
image of $m$ points. Let $(Y,\lambda)$
be a digital interval or a digital
simple closed curve of $n$ points, with $n>m+2$.
Then there is no universal function from $X$ to $Y$.
\end{Proposition}

{\em Proof}: Let $f: X \rightarrow Y$ be a
$(\kappa,\lambda)$ continuous function. Then
$f(X)$ is a $\lambda$-connected subset of $Y$, and
$|f(x)| \leq m < n$.

We show that $Y \setminus f(X)$ has a component with at
      least 2 points, one of which is not 
      $\lambda$-adjacent to any member of $f(X)$.
\begin{itemize}
\item If $Y$ is a digital interval $[a,a+n-1]_{\bf Z}$, then,
      since $f(X)$ is a connected subset of $Y$,
      $f(X)=[u,v]_{\bf Z}$. Consider the
      following possibilities.
      \begin{itemize}
      \item $v \leq a+n-3$. Then the endpoint 
            $a+n-1$ of $Y\setminus f(X)$ is
             not adjacent to any point of $f(X)$.
      \item $v > a+n-3$. Therefore, $v \geq a+n-2$. Then
            \[ u=v-|f(X)|+1 \geq a+n-2-|f(X)|+1 \geq a+n-(m+2)+1
               > a+1.
            \]
            I.e., $u \geq a+2$, so the point $a$ of 
            $Y\setminus f(X)$ is not adjacent to any 
            point of $f(X)$.
      \end{itemize}
      \item If $Y$ is a digital simple closed curve, we may
            assume $Y = \{y_j\}_{j=0}^{n-1}$, where 
            $y_a$ and $y_b$ are adjacent if and only if
            $a=(b+1)~mod~n$ or $a=(b-1)~mod~n$. Since
            $f(X)$ is connected, we may assume without loss of
            generality that $f(X)=\{y_j\}_{j=0}^r$ where
            $0 \leq r < m < n-2$. Then $y_{r+2}$ is a point of
            $Y \setminus f(X)$ that is not adjacent to any point of $f(X)$.
\end{itemize}
Thus, $f(X)$ is not $\lambda$-dominating in $Y$. The
assertion follows from Theorem~\ref{univ-dominating}.
$\Box$

\begin{Proposition}
\label{AFPP-id-univ}
Let $(X,\kappa)$ be a digital image. Then $(X,\kappa)$
has the AFPP if and only if the identity function 
$1_X$ is universal for $(X,X)$.
\end{Proposition}

{\em Proof}: The function $1_X$ is universal if and only if
for every
$(\kappa,\kappa)$-continuous $f: X \rightarrow X$,
there exists $x \in X$ such that 
$f(x) \leftrightarrow_{\kappa} 1_X(x)=x$, which is true
if and only if $(X,\kappa)$ has the AFPP.
$\Box$

\begin{Theorem}
\label{restriction}
Let $(X,\kappa)$ and $(Y,\lambda)$ be digital images
and let $U \subset X$. If the restriction function
$f|_U \, : (U,\kappa) \rightarrow (Y,\lambda)$ is a
universal function for $(U,Y)$, then $f$ is a
universal function for $(X,Y)$.
\end{Theorem}

{\em Proof}: Let $h: X \rightarrow Y$ be
$(\kappa,\lambda)$-continuous. Since $f|_U$ is
universal, there exists $u \in U \subset X$ such that
$h(u)=h|_U(u) \leftrightarrow_{\lambda} f|_U(u)=f(u)$. 
Hence $f$ is universal for $(X,Y)$. $\Box$

\begin{Theorem}
\label{univ-composition}
Let $(W,\kappa)$, $(X,\lambda)$, and $(Y,\mu)$ be
digital images. Let
$f: W \rightarrow X$ be $(\kappa,\lambda)$-continuous
and let
$g: X \rightarrow Y$ be $(\lambda,\mu)$-continuous.
If $g \circ f$ is universal, then $g$ is also
universal.
\end{Theorem}

{\em Proof}: Let $h: X \rightarrow Y$ be 
$(\lambda,\mu)$-continuous. Since $g \circ f$ is
universal, there exists $w \in W$ such that
$(g \circ f)(w) \leftrightarrow_{\mu} (h \circ f)(w)$. I.e., for 
$x=f(w) \in X$ we have
$g(x) \leftrightarrow_{\mu} h(x)$.
Since $h$ was arbitrarily chosen, the assertion
follows. $\Box$

\begin{Theorem}
\label{univ-equivalents}
If $g: (U,\mu) \rightarrow (X,\kappa)$ and
$h: (Y, \lambda) \rightarrow (V,\nu)$ are 
digital isomorphisms and $f: X \rightarrow Y$ is
$(\kappa,\lambda)$-continuous, then the following
are equivalent.
\begin{enumerate}
\item $f$ is a universal function for $(X,Y)$.
\item $f \circ g$ is universal.
\item $h \circ f$ is universal.
\end{enumerate}
\end{Theorem}

{\em Proof}: (1 implies 2): Let
$k: U \rightarrow Y$ be $(\mu,\lambda)$-continuous.
Since $f$ is universal, there exists $x \in X$ such
that $(k \circ g^{-1})(x) \leftrightarrow f(x)$.
By substituting $x=g(g^{-1}(x))$,
we have $k(g^{-1}(x)) \leftrightarrow (f \circ g)(g^{-1}(x))$. 
Since $k$ was
arbitrarily chosen and $g^{-1}(x) \in U$, it follows
that $f \circ g$ is universal.

(2 implies 1): This follows from
Theorem~\ref{univ-composition}.

(1 implies 3): Let $m: X \rightarrow V$ be
$(\kappa,\nu)$-continuous. Since $f$ is universal,
there exists $x \in X$ such that 
$(h^{-1} \circ m)(x) \leftrightarrow f(x)$.
Then $m(x)=h((h^{-1} \circ m)(x)) \leftrightarrow_{\nu}
(h \circ f)(x)$. Since $m$ was
arbitrarily chosen, it follows that $h \circ f$ is
universal.

(3 implies 1): Suppose $h \circ f$ is universal.
Then given a $(\kappa, \lambda)$-continuous
$r: X \rightarrow Y$, there exists $x \in X$ such
that $h \circ f(x) \leftrightarrow_{\nu} h \circ r(x)$. 
Therefore,
$f(x)=(h^{-1} \circ h \circ f)(x) \leftrightarrow_{\lambda}
(h^{-1} \circ h \circ r)(x)=r(x)$. Since $r$ was arbitrarily
chosen, it follows that $f$ is universal.
$\Box$

\begin{Corollary}
\label{iso-equivs}
Let $f: (X,\kappa) \rightarrow (Y,\lambda)$ be a
digital isomorphism. Then $f$ is universal for
$(X,Y)$ if and only if $(X,\kappa)$ has the AFPP.
\end{Corollary}

{\em Proof}: The function $f$ is universal, by
Theorem~\ref{univ-equivalents}, if and only if
$f \circ f^{-1}=1_X$ is universal, which, by
Proposition~\ref{AFPP-id-univ}, is true if and
only if $(X,\kappa)$ has the AFPP. $\Box$

It may be useful to remind the reader for the following
theorem that points that are $c_n$-adjacent in ${\bf Z}^n$
may differ in every coordinate.
Concerning products, we have the following.

\begin{Theorem}
\label{products}
Let $(X_i,c_{n_i}) \subset {\bf Z}^{n_i}$,
$i=1,2, \ldots, m$. Let $s=\sum_{i=1}^m n_i$.
Consider the digital image
$X = \Pi_{i=1}^m \, X_i \subset {\bf Z}^s$. If
$(X,c_s)$ has the AFPP then each $(X_i,c_{n_i})$
has the AFPP.
\end{Theorem}

{\em Proof}: Suppose $(X,c_s)$ has the AFPP.
Let $f_i: X_i \rightarrow X_i$ be
$(c_{n_i}, c_{n_i})$-continuous. Then the function
$f: X \rightarrow X$ defined by
\[ f(x_1, x_2, \ldots, x_m)=
   (f_1(x_1), f_2(x_2), \ldots, f_m(x_m))
\]
is $(c_s,c_s)$-continuous. By
Proposition~\ref{AFPP-id-univ}, $1_X$ is universal
for $(X,X)$. Therefore, there is a point
$x_* = (x_{1,*}, x_{2,*}, \ldots, x_{m,*}) \in X$
with $x_{i,*} \in X_i$ such that 
$x_* \leftrightarrow_{c_s} f(x_*)$. Therefore,
$x_{i,*} \leftrightarrow_{c_{n_i}} f_i(x_{i,*})$ 
for all $i$. Since $f_i$ was
arbitrarily chosen, it follows that $1_{X_i}$ is
universal for $(X_i,X_i)$.
Proposition~\ref{AFPP-id-univ} therefore implies
that $(X_i,c_{n_i})$ has the AFPP. $\Box$

\section{Summary}
\label{summary-sec}
We have shown that only single-point digital images
have the fixed point property. However, digital 
$n$-cubes have the approximate fixed point property
with respect to the $c_n$-adjacency
~\cite{Rosenfeld}.
We have shown that the approximate fixed point 
property is preserved
by digital isomorphism and by digital retraction, and we
have a result concerning preservation of the AFPP by
Cartesian products. We
have studied relations between universal functions
and the AFPP. We have corrected several errors that
appeared in previous papers.

\section{Acknowledgment}
\label{ack-sec}
The remarks of an anonymous reviewer were very helpful and are
gratefully acknowledged.


\begin{thebibliography}{XX}
\bibitem{AKO}
H. Arslan, I. Karaca, and A. \u{O}ztel,
Homology groups of $n$-dimensional digital images,
XXI {\em Turkish National Mathematics Symposium} (2008),
B1-13.

\bibitem{Borsuk}
K. Borsuk,
{\em Theory of Retracts},
Polish Scientific Publishers, Warsaw, 1967.

\bibitem{Boxer}
L. Boxer,
Digitally continuous functions,
{\em Pattern Recognition Letters} 15 (1994), pp. 833-839, 1994.

\bibitem{Boxer99}
L. Boxer,
A classical construction for the digital fundamental group,
{\em Journal of Mathematical Imaging and Vision} 10 (1999), pp. 51-62, 1999.

\bibitem{Boxer05}
L. Boxer,
Properties of digital homotopy,
{\em Journal of Mathematical Imaging and Vision} 22 (2005), 19-26.

\bibitem{Boxer-prod}
L. Boxer,
Digital products, wedges, and covering spaces,
{\em Journal of Mathematical Imaging and Vision} 25 (2006), 159-171.

\bibitem{Boxer07}
L. Boxer,
Fundamental Groups of Unbounded Digital Images,
{\em Journal of Mathematical Imaging and Vision}
27 (2007), 121-127. 

\bibitem{Boxer10}
L. Boxer,
Continuous Maps on Digital Simple Closed Curves,
{\em Applied Mathematics} 1 (2010), 377-386. 

\bibitem{BKO}
L. Boxer, I. Karaca, and A \u{O}ztel,
Topological invariants in digital images,
{\em Journal of Mathematical Sciences: Advances and
Applications} 11 (2) (2011), 109-140.

\bibitem{C&L}
G. Chartrand and L. Lesniak,
{\em Graphs $\&$ Digraphs}, 2nd ed.,
Wadsworth, Inc., Belmont, CA, 1986.

\bibitem{Chen1}
L. Chen,
Gradually varied surfaces and its
optimal uniform approximation,
{\em SPIE Proceedings} 2182 (1994),
300-307

\bibitem{Chen2}
L. Chen,
{\em Discrete Surfaces and Manifolds},
Scientific Practical Computing,
Rockville, MD, 2004

\bibitem{D&K}
E. Demir and I. Karaca,
Simplicial homology groups of certain digital surfaces, 
{\em Hacettepe Journal of Mathematics and Statistics}, 
to appear.

\bibitem{Ege&Karaca}
O. Ege and I. Karaca,
Lefschetz fixed point theorem for digital images,
{\em Fixed Point Theory and Applications} 2013, 2013:253.
Available at $\\$
http://www.fixedpointtheoryandapplications.com/content/2013/1/253

\bibitem{Ege&Karaca2}
O. Ege and I. Karaca,
Fundamental properties of digital simplicial homology groups,
{\em American Journal of Computer Technology and Application}
1 (2) (2013), 25-42.

\bibitem{Ege&Karaca3}
O. Ege and I. Karaca,
Applications of the Lefschetz Number to Digital Images,
{\em Bulletin of the Belgian Mathematical Society, Simon Stevin}
21 (5) (2014), 823-839.

\bibitem{Han}
S.-E. Han,
Non-product property of the digital fundamental group,
{\em Information Sciences} 171 (2005), 73-91.

\bibitem{Han07}
S.-E. Han,
Digital fundamental group and Euler characteristic of a connected
sum of digital closed surfaces,
{\em Information Sciences} 177 (16) (2007), 3314-3326.

\bibitem{Herman}
G.T. Herman,
Oriented surfaces in digital spaces,
{\em CVGIP: Graphical Models and Image Processing} 55, pp. 381-396, 1993.

\bibitem{Karaca&Ege}
I. Karaca and \u{O}. Ege,
Some results on simplicial homology groups of 2D digital
images, {\em International Journal of Information and
Computer Science} 1 (8) (2012), 198-203.

\bibitem{Khalimsky}
E. Khalimsky,
{\em Motion, deformation, and homotopy in finite spaces}, in
Proceedings IEEE Intl. Conf. on Systems, Man, and Cybernetics,
pp. 227-234, 1987.

\bibitem{Kong}
T.Y. Kong,
A digital fundamental group,
{\em Computers and Graphics} 13, pp. 159-166, 1989.



\bibitem{Rosenfeld}
A. Rosenfeld,
`Continuous' functions on digital pictures,
{\em Pattern Recognition Letters} 4, pp. 177-184, 1986.

\bibitem{Spanier}
E.H. Spanier,
{\em Algebraic Topology},
McGraw-Hill, New York, 1966.

\end{thebibliography}
\end{document}